\documentclass[a4paper,12pt]{amsart}

\usepackage{amsmath}
\usepackage{amsfonts}
\usepackage{amssymb}
\usepackage{graphicx}
\usepackage{hyperref}
\usepackage{cite}

\usepackage{amsthm}
\usepackage[all,cmtip]{xy}

\newcommand{\mbQ}{\mathbb{Q}}

\newcommand{\mbP}{\mathbb{P}}

\newcommand{\mcO}{\mathcal{O}}

\DeclareMathOperator{\Aut}{Aut}

\DeclareMathOperator{\Mov}{Mov}
\DeclareMathOperator{\lct}{lct}
\DeclareMathOperator{\Nef}{Nef}

\DeclareMathOperator{\Pic}{Pic}
\DeclareMathOperator{\Bir}{Bir}

\theoremstyle{plain}
\newtheorem{theorem}{Theorem}[section]

\newtheorem{question}[theorem]{Question}

\theoremstyle{definition}
\newtheorem{definition}[theorem]{Definition}

\theoremstyle{remark}
\newtheorem{remark}[theorem]{Remark}

\title[Rationally connected varieties and Mori dream spaces]{A remark on rationally connected varieties and Mori dream spaces}

\author{Claudio Fontanari}
\address{Claudio Fontanari \\ Dipartimento di Matematica, Universit\`a di Tren\-to, Via Sommarive 14, 38123 Povo, Trento, Italy.}
\email{claudio.fontanari@unitn.it}

\author{Diletta Martinelli}
\address{Diletta Martinelli\\ School of Mathematics, University of Edinburgh, James Clerk Maxwell Building, The King's Buildings, Peter Guthrie Tait Road, Edinburgh, EH9 3FD.}
\email{Diletta.Martinelli@ed.ac.uk}

\begin{document}
	\bibliographystyle{alpha}

\begin{abstract}
	In this short note, we show that a construction by Ottem \cite[Theorem 1.1]{Ott15} provides an example of a rationally connected variety that is not birationally equivalent 
	to a Mori dream space with terminal singularities. This answers in the negative (at least in the category of terminal varieties) a question posed by Krylov \cite[Remark 5.7]{Kry14}.  
\end{abstract}

\subjclass[2010]{Primary 14E30; Secondary 14E07}
\keywords{Mori dream space, rationally connected variety, birational rigidity}

\maketitle


\section{Introduction}

Varieties of Fano type are examples of varieties that behave well with respect to the Minimal Model Program. They are known to be rationally connected by \cite{KMM92} and \cite{Zha06}. However, the converse is not true (the blow-up of $\mbP^2$ in 10 very general points provides an obvious counterexample). In the recent paper \cite{Kry14} it is shown that there exist smooth rationally connected varieties of dimension $n \ge 4$ that are not birationally equivalent to a variety of Fano type.

Mori dream spaces form another class of varieties that behave well with respect to a $D$-MMP, for any divisor $D$ \cite[Proposition 1.11]{HK00}. 
We recall (see \cite[Definition 1.10]{HK00}) that a Mori dream space is a normal $\mathbb{Q}$-factorial projective variety such that
\begin{itemize}
	\item[(i)] $\Pic(X)$ is finitely generated;
	\item[(ii)] The Nef cone $\Nef(X)$ is the affine hull of finitely many semiample line bundles;
	\item[(iii)] There is a finite dimensional collection of small $\mathbb{Q}$-factorial mo\-difications $f_i \colon X \dashrightarrow X_i$ such that each $X_i$ satisfies (ii) and the movable cone $\Mov(X)$ is the union of the $f^*(\Nef(X_i))$.
\end{itemize}
 
 It was proven in \cite[Corollary 1.3.2]{BCHM10} that any $\mbQ$-factorial projective variety of Fano type is a Mori dream space.
 Krylov then asked the following question.
\begin{question}\cite[Remark 5.7]{Kry14} \label{quest:Kry}
	Let $X$ be a rationally connected variety. Is $X$ birationally equivalent to a Mori dream space?
\end{question} 
 
In this short note, we claim that a negative answer to Question \ref{quest:Kry} is implied 
(at least in the category of terminal varieties)
by \cite[Theorem 1.1]{Ott15}, stating that a very general hypersurface of bidegree $(d, e)$ in $\mbP^1 \times \mbP^n$ is not a Mori dream space for $d \ge n+1$ and $e \ge 2$.

More precisely, we prove the following fact.
\begin{theorem} \label{thm:main}
For every $n \ge 11$ and $d \ge n+1$ there exists a smooth very general hypersurface $X$ in $\mbP^1 \times \mbP^n$ of bidegree $(d, n)$ which is rationally connected but not birationally equivalent to a Mori dream space with terminal singularities.
\end{theorem}

In Section \ref{sec:preli} we recall the necessary notions from the Minimal Model Program and the definition of birationally rigid varieties. In Section \ref{sec:proof} we prove Theorem \ref*{thm:main}.
The strategy of the proof is quite simple: since we start from a variety $X$ that is not a Mori dream space, we only need to ensure that $X$ is birationally superrigid and it does not admit fibre-wise transformations.

\section{Preliminaries} \label{sec:preli}
Throughout the paper we work over the field of complex numbers. All the varieties we consider are assumed to be normal projective and $\mbQ$-factorial.

\subsection{Minimal Model Program}
We recall the standard definition of singularities appearing in the Minimal Model Program. For more details see \cite[Section 2.3]{KM08}.

\begin{definition}\cite[Definition 2.34]{KM08}
	Let $X$ be a normal variety and let $\Delta$ be an effective $\mbQ$-divisor on $X$. Let $\pi \colon \tilde{X} \to X$ be a birational morphism from a normal variety $\tilde{X}$. Let $\tilde{\Delta} = \pi^{-1}_*(\Delta)$ be the proper transform of $\Delta$.
	Then we can write 
	\begin{equation*}
		K_{\tilde{X}} + \tilde{\Delta} = \pi^*(K_X + \Delta) + \sum\limits_E a(E, X, \Delta)
	\end{equation*}
	where $E$ runs through all the distinct exceptional prime divisors on $\tilde{X}$ and $a(E, X, \Delta)$ is a rational number. We say that the pair $(X, \Delta)$ is terminal (resp. canonical, log terminal, log canonical) if $a(E, X, \Delta) > 0$ (resp. $a(E, X, \Delta) \geq 0$, $a(E, X, \Delta) > -1$, $a(E, X, \Delta) \geq -1$) for every prime divisor $E$ on $\tilde{X}$. If $\Delta = 0$ then we simply say that $X$ has terminal (resp. canonical, log terminal, log canonical) singularities.
\end{definition}

We now define the log canonical threshold of a pair (see for details \cite[Section 8]{Kol97}).

\begin{definition} \cite[Definition 1.2]{Che09}
	Let $X$ be a variety with at most log terminal singularities, let $Z \subseteq X$ be a closed subvariety, and let $D$ be an effective $\mbQ$-Cartier $\mbQ$-divisor on $X$. Then the number 
	\begin{equation*}
	\lct_Z(X, D) = \sup \{\lambda \in \mbQ | \text{ the log pair }(X, \lambda D) \text{ is log canonical along }Z\} 
	\end{equation*} 
	is said to be the log canonical threshold of $D$ along $Z$. We assume, in addition, that $X$ is a Fano variety. We then define the log canonical threshold of $X$ by the number
$$
		\lct(X) = \inf \{ \lct(X, D) | D \text{ is an effective } \mbQ\text{-divisor on } X \text{ s.t. } D \equiv -K_X \}.
$$ 
\end{definition}

The number $\lct(X)$ is an algebraic counterpart of the so-called $\alpha$-invariant first introduced by Tian in \cite{Tia87}.  

\subsection{Birational rigidity}

\begin{definition} \label{def:mfs}
	A Mori fiber space is a $\mbQ$-factorial projective variety $X$ with at most terminal singularities and a morphism $\phi \colon X \to Z$, such that
	\begin{itemize}
		\item The anticanonical class of $X$, $-K_X$, is $\phi$-ample;
		\item The relative Picard number, $\Pic(X/Z)$, is 1;
		\item $\dim Z < \dim X$.
	\end{itemize}
\end{definition}

Fano varieties with Picard rank 1 and Fano fibrations over $\mbP^1$, by which we mean terminal $\mbQ$-factorial varieties with Picard number 2 and a map to $\mbP^1$ such that the generic fiber is a smooth Fano variety, are typical examples of Mori fiber spaces.

We recall here just the definition of birationally superrigidity, while for a comprehensive introduction to the subject we refer to \cite{Puk13} and \cite{Che05}.

\begin{definition} \cite[Definition 1.3]{CM04}
	Let $X \to Z$ and $X' \to Z'$ two Mori fiber spaces, a birational map $f \colon X \dashrightarrow X'$ is \emph{square} if fits into the commutative diagram
\begin{equation*}
	\xymatrix{
		X \ar[d] \ar@{-->}[r]^f & X'\ar[d]\\
		Z \ar@{-->}[r]^g        & Z'}
\end{equation*}	
where $g$ is birational and the map induced on the generic fiber $f_L \colon X_L \to Z_L$ is biregular, where we denote with $L$ the generic point of $Z$. In this case we say that $X/Z$ and $X'/Z'$ are \emph{square equivalent}. 
\end{definition}
\begin{definition}
We say that a Mori fiber space is \emph{birationally rigid} if the set 
\begin{equation*}
\{\text{Mori fiber space } Y \to S | Y \text{ birational to } X\}/\text{ square equivalence}
\end{equation*}
contains just a single element. Moreover, we say that $X$ is \emph{birationally superrigid} if in addition the group of birational automorphisms $\Bir(X)$ and the group of biregular automorphisms $\Aut(X)$ coincide.
\end{definition}

Therefore, it follows that if $X/Z$ and $X'/Z'$ are Mori fiber spaces and $f \colon X \dashrightarrow X'$ is a birational map between them, then $f$ maps $X$ to $X'$ fibre-wise.

\section{Proof of Theorem \ref{thm:main}} \label{sec:proof}

\begin{remark}
	The hypersurface $X$ admits a fibration onto  $\mbP^1$, whose generic fiber is a Fano variety by the adjunction formula. 
	Hence $X$ is rationally connected by \cite[Corollary 1.3]{GHS03}.
\end{remark}

Let $U \subset \mbP H^0(\mbP^1 \times \mbP^n, \mcO(d,n))$ be the dense set corresponding to hypersurfaces $f$ 
which are not Mori dream spaces by \cite{Ott15}.  

On the other hand, by \cite[Theorem 4]{Puk15}, if $n \ge 11$ then there exists a Zariski open subset $\mathcal{F}_{\mathrm{reg}} \subset 
\mbP H^0(\mbP^n, \mcO(n))$ with complement of codimension $> 1$ such that every hypersurface $F \in \mathcal{F}_{\mathrm{reg}}$ satisfies:

(i) $F$ is a factorial Fano variety with terminal singularities and $\mathrm{Pic}(F)$ $= \mathbb{Z} K_F$;

(ii) for every effective divisor $D \in \vert -K_F \vert$ the pair $(F,\frac{1}{n}D)$ is log canonical, and for every mobile linear system 
$\Sigma \subset \vert -K_F \vert$ the pair $(F,\frac{1}{n}D)$ is canonical for a general divisor $D \in \Sigma$. 

In particular, this means that $\lct(F) \geq 1$.
\\

We consider the natural evaluation and projection maps: 
\begin{align*}
	ev: \mbP H^0(\mbP^1 \times \mbP^n, \mcO(d,n)) \times \mbP^1 &\to \mbP H^0(\mbP^n, \mcO(n))\\
	(f,p) &\mapsto f(p)\\
	\pi: \mbP H^0(\mbP^1 \times \mbP^n, \mcO(d,n)) \times \mbP^1 &\to \mbP H^0(\mbP^1 \times \mbP^n, \mcO(d,n))\\ 
	(f,p) &\mapsto f
\end{align*}
and let 
$$
V := \mbP H^0(\mbP^1 \times \mbP^n, \mcO(d,n)) \setminus \pi(ev^{-1}(\mbP H^0(\mbP^n, \mcO(n)) \setminus \mathcal{F}_{\mathrm{reg}})).
$$
The set $V$ is a Zariski open subset since $\mathcal{F}_{\mathrm{reg}}$ is so and it is non-empty since the complement of $\mathcal{F}_{\mathrm{reg}}$ has codimension $>1$.

Now, if $f \in U \cap V \ne \emptyset$ then the Mori fiber space $X$ defined by $f$ is birationally superrigid (see for instance \cite[Proposition 3.1, pp. 309--310]{Puk13}: 
as in \cite[Lemma 3.7]{Kry14}, the K-condition is trivially satisfied for $d >>0$). We can also exclude 
fibre-wise tranformations by quoting \cite[Theorem 1.5]{Che09}, exactly as in \cite[Corollary 3.2]{Kry14}. 
It follows that $X$ is not birational to a Mori dream space with terminal singularities. Indeed, if $Y$ were a Mori dream space birational 
to $X$, then since $X$ has negative Kodaira dimension $Y$ would be birational via a Minimal Model Program to a Mori fiber space preserving the structure of Mori dream space, a contradiction. 

\subsection{Open Questions}
If we start from a rationally connected variety and we run a MMP, we end up with a Mori fiber space as in Definition \ref{def:mfs}. Therefore, an interesting question related to the previous results is the following.

\begin{question}
	Which Mori fiber spaces over $\mbP^1$ are Mori dream spaces? Is it possible to reach some kind of classification?
\end{question}  
In dimension two, Mori fiber spaces over $\mbP^1$ are the Hirzebruch surfaces, that are toric and, therefore, Mori dream spaces.
\\

Further connections between Mori Dream Spaces and the birational geometry of Fano varieties are suggested in \cite{AZ16}.

\section*{Acknowledgements}
We would like to thank Ivan Cheltsov and John Ottem for useful conversations on this subject.
The first named author is partially supported by GNSAGA of INdAM, by PRIN 2015 
"Geometria delle variet\`a algebriche", and by FIRB 2012 "Moduli spaces and Applications". The second named author was supported by the ERC starting grant WallXBirGeom 337039.

\end{document}